\documentclass[12pt]{article}
\usepackage{latexsym}
\usepackage{amsmath, amscd, amssymb}
\usepackage{eepic,epic}
\begin{document}
\baselineskip 16pt

\def \Aut {\mbox{\rm Aut\,}}
\def \Rot {\mbox{\rm Rot\,}}
\def \Inv {\mbox{\rm Inv\,}}
\def \Mon {\mbox{\rm Mon\,}}
\def \Ex {\mbox{\rm Ex\,}}
\def \Exo {\mbox{\rm Exo\,}}
\newfont{\sBbb}{msbm7 scaled\magstephalf}
\newcommand{\BC}{\mbox{\Bbb C}}
\newcommand{\BR}{\mbox{\Bbb R}}
\newcommand{\BZ}{\mbox{\Bbb Z}}
\newcommand{\sZ}{\mbox{\sBbb Z}}

\newcommand{\Iso}{\mbox{Iso}\,}
\newcommand{\Fi}{\mbox{Fix}\,}
\renewcommand{\labelenumi}{\theenumi}
\newcommand{\qed}{\mbox{\raisebox{0.7ex}{\fbox{}}}}
\newtheorem{theorem}{Theorem}[section]
\newtheorem{problem}{Problem}
\newtheorem{defin}{Definition}
\newtheorem{lemma}[theorem]{Lemma}
\newtheorem{prop}[theorem]{Proposition}
\newtheorem{conj}{Conjecture}
\newtheorem{op}{Open Problem}
\newtheorem{example}{Example}[section]
\newtheorem{note}{Note}
\newtheorem{remark}{Remark}
\newtheorem{corollary}[theorem]{Corollary}
\newenvironment{pf}{\medskip\noindent{Proof:}
  \hspace{-.5cm}      \enspace}{\hfill \qed \newline \smallskip}
\newenvironment{pflike}[1]{\medskip\noindent {#1}
   \enspace}{\medskip}
\newenvironment{defn}{\begin{defin} \em}{\end{defin}}
\newenvironment{nt}{\begin{note} \em}{\end{note}}
\newenvironment{rem}{\begin{remark} \em}{\end{remark}}
\newenvironment{examp}{\begin{example} \rm}{\end{example}}
\newcommand{\lra}{\longrightarrow}
\newcommand{\vect}[2]{\mbox{$({#1}_1,\ldots,{#1}_{#2})$}}
\newcommand{\comb}[2]{\mbox{$ \left( \begin{array}{c}
        {#1} \\ {#2} \end{array}\right)$}}
\newcommand{\tve}[1]{\mbox{$\mathbf{\tilde {#1}}$}}
\newcommand{\bve}[1]{\mbox{$\mathbf{{#1}}$}}
\setlength{\unitlength}{12pt}
\renewcommand{\labelenumi}{(\theenumi)}
\def\mod{\hbox{\rm mod }}

\title{
Classification of  nonorientable regular embeddings of
Hamming graphs}
\author{Gareth A. Jones \\
 {\small School of Mathematics,
University of Southampton, Southampton SO17 1BJ, U.K.}  \\ 
 Young Soo Kwon     \\
{\small Mathematics, Pohang University of Science and Technology,
Pohang, 790-784 Korea}\vspace{5mm}
 }

\date{}
\maketitle

\renewcommand{\thefootnote}{*}

\begin{abstract}
By a regular embedding of a graph $K$ in a surface we mean a
2-cell embedding of $K$ in a compact connected surface such that the
automorphism group acts regularly on flags. In this paper, we classify the nonorientable regular embeddings of the Hamming graph $H(d,n)$.  We show that there exists such an embedding if and only if $n=2$ and $d=2$, or $n= 3$ or $4$ and $d\ge 1$, or $n=6$ and $d=1$ or $2$. We also give constructions and descriptions of these embeddings.
 \vskip10pt

\end{abstract}

\section{Introduction}

 A {\it map} $\cal M$ is a 2-cell embedding of a
graph $K$ in a compact, connected surface $S$. An {\it automorphism} of $\cal M$
is a permutation of its {\em flags\/} (mutually
incident vertex-edge-face triples) which preserves their relations of having a vertex, edge or face in common; it therefore induces an automorphism of $K$ which extends to a
self-homeomorphism of $S$.
The group $G:=\Aut({\cal M})$ of all automorphisms of $\cal M$ acts semi-regularly on its flags, so $|G|\le 4|E|$, where $E$ is the set of edges. If this bound is attained then $G$ acts
regularly on the flags, and $\cal M$ is called a {\it regular} map.
Equivalently, $\cal M$ is regular if and only if there are three
involutions $\lambda$, $\rho$ and $\tau$ in
$G$, each fixing a distinct pair of elements $v,e,f$ of some
flag $(v,e,f)$; in this case we have $G =\langle
\lambda,\rho,\tau\rangle$. In what follows we shall assume that $\cal M$ is a regular map, with
$\lambda$ fixing $e$ and $f$, and $\rho$ fixing $v$ and $f$, so $\tau$ fixes $v$ and $e$. We call such a triple $(\lambda, \rho, \tau)$ an \emph{admissible triple}.

Since $\tau\lambda=\lambda\tau$ the stabilizer $G_e=\langle \lambda, \tau \rangle$ in $G$ of $e$ is a dihedral group of order 4, i.e.~a Klein four-group, isomorphic to ${\mathbb Z}_2\times
{\mathbb Z}_2$. Similarly, the stabilizers $G_v=\langle \rho, \tau \rangle$ and $G_f=\langle \lambda, \rho \rangle$ of $v$ and $f$ are dihedral groups of orders
$2m$ and $2n$, where $m$ is the common valency of the vertices of $\cal M$, i.e.~the order of $\rho\tau$, and $n$ is the covalency (the number of edges of each face), equal to the order of $\lambda\rho$.

When a regular map $\cal M$ is represented in this way by a triple of involutions
$\lambda,\rho,\tau$ we write ${\cal M}={\cal
M}(\lambda,\rho,\tau)$. Two regular maps ${\cal
M}(\lambda,\rho,\tau)$ and ${\cal M'}(\lambda',\rho',\tau')$ with
underlying graphs $K$ and $K'$ are {\em isomorphic} if there is a
graph isomorphism  $\psi:K\to K'$ such that $\psi
^{-1}\lambda\psi= \lambda'$, $\psi^{-1}\rho\psi=\rho'$ and
$\psi^{-1}\tau\psi=\tau'$. A detailed explanation of the above
representation of regular maps can be found in \cite[Theorem
3]{GNSS}. The basic theory of regular maps, as well as other relevant
information, can also be found in~\cite{JS1,JoTh,N,Vi1,Vi2}.

 If a regular map $\cal M$ is obtained from an embedding $i:K\to S$ of a graph $K$
in a surface $S$ we say that $i$ is a {\it regular embedding\/} of
$K$. The surface $S$ underlying $\cal M$ is nonorientable if and
only if there is a cycle $C$ in $K$ with a neighbourhood in $S$
homeomorphic to a M\" obius band. Such a cycle will be
called {\it orientation-reversing}. A regular map $\cal M$ is
nonorientable if and only if its automorphism group $G$
is generated by $R := \rho\tau$ and the involution $L := \lambda\tau=\tau\lambda$. In particular, if
there is an orientation-reversing cycle $C$ of length $l$ in
$\cal M$ then there is an associated relation of the form
$LR^{m_1}LR^{m_2}\dots LR^{m_l}=\tau$ in $G$.
Conversely, the existence of such a relation in $G$
implies that $\cal M$ is nonorientable.  We call such a triple $(\lambda, \rho, \tau)$ an \emph{nonorientable admissible triple}.

There are only a few families of graphs for which a complete classification of
their nonorientable regular embeddings is known. Such embeddings of complete graphs $K_n$ have been classified by
James~\cite{Ja} and Wilson~\cite{Wi}: these exist if and only if $n$ is
$3, 4$ or $6$.  Nedela and the second
author~\cite{KN} have shown the nonexistence of a nonorientable
regular embedding of the $n$-dimensional cube graph $Q_n$ for all $n$ except  $n =2$.  In contrast with all other known
cases, the complete bipartite graph $K_{n,n}$ has a nonorientable regular embedding
for infinitely many values of $n$, as shown by Kwak and the second author~\cite{KK1}: in fact, such an embedding
exists if and only if $n \equiv 2$ mod $(4)$ and all odd prime
divisors of $n$ are congruent to $\pm 1$ mod $(8)$.

A map $\cal M$ is {\em orientably regular} if the underlying surface is orientable and the orientation-preserving subgroup $\Aut^+({\cal M})$ of $\Aut({\cal M})$ acts regularly on the arcs (directed edges) of $\cal M$. The first author~\cite{J} has recently obtained a classification of such embeddings of Hamming graphs $H(d,n)$, and this includes a classification of their orientable regular embeddings (called {\em reflexible\/} embeddings there).  In this paper, we aim to classify the nonorientable regular embeddings of Hamming graphs.  Our main
result is the following theorem:

\begin{theorem} \label{main-theorem}
There exists a nonorientable regular embedding of the Hamming graph
$H(d,n)$ if and only if either
\begin{enumerate}
\item $n=2$ and $d=2$, or
\item  $n=3$ or $4$ and $d\ge 1$, or
\item $n=6$ and $d=1$ or $2$.
\end{enumerate}
In cases~(1) and (2) the embedding of $H(d,n)$ is unique up to isomorphism, whereas in case~(3) there are two such embeddings for each of the two graphs $H(d,n)$.
\end{theorem}

Further information about each of these maps, namely its type, genus and automorphism group, is given in Section~2.

This paper is organized as follows. In Section~2 we construct and describe some examples of nonorientable regular embeddings of $H(d,n)$, and in Section~3 we classify all such embeddings by showing that each of them is isomorphic to one of these examples.

The authors are grateful to the organisers of GEMS 09 in T\'ale, Slovakia, and of the Algebraic Graph Theory Summer School in Rogla, Slovenia, 2011, for providing the opportunities for this collaboration.

\section{Construction of nonorientable Hamming maps}

The \emph{Hamming graph\/} $H(d,n)$ is the cartesian product of $d$ cliques of size $n$. Specifically, we define $H(d,n)$ to have vertex set $V=[n]^d$ where $[n]=\{0, 1, \ldots. n-1\}$ for some $n \ge 2$, with two vertices $u=(u_i )$ and $v=(v_i)$ adjacent if and only if $u_i = v_i$ for all except exactly one value of $i$.  

The automorphism group $\Aut(H(d,n))$ of this graph is the wreath
product $S_n \wr S_d$ of the symmetric groups $S_n$ and $S_d$.
This is a semidirect product of a normal subgroup $S_n \times S_n
\times \cdots \times S_n$, whose $i$-th direct factor acts on
$i$-th coordinate of each vertex and fixes the $j$-th coordinates for
$j \neq i$, by a complement $S_d$ which permutes the $d$ coordinates of each vertex.

Following Coxeter and Moser~\cite[Ch.~8]{CM}, we say that a regular map has {\em type} $\{p, q\}_r$ if $p$ is its covalency (the number of sides of each face), $q$ is the valency of each vertex, and $r$ is its Petrie length (the length of each Petrie polygon).

In the case $d=1$ the regular embeddings of $H(d,n)$ are already known: $H(1,n)$ is isomorphic to the complete graph $K_n$, and the regular embeddings of this graph have been classified by James and the first author~\cite{JJ} in the orientable case, and by James~\cite{Ja} and Wilson~\cite{Wi} in the nonorientable case. The results are as follows:

\begin{prop}\label{regKn}
\noindent {\rm (a)} There are, up to isomorphism, just three orientable regular embeddings of complete graphs $K_n$ for $n\ge 2$ ; they are embeddings of a closed interval, a triangle and a tetrahedron in the sphere, with $n=2, 3$ and $4$.

\smallskip

\noindent {\rm (b)} There are, up to isomorphism, just four nonorientable regular embeddings of complete graphs $K_n$; they are the antipodal quotients of a hexagon and a cube on the sphere, giving regular embeddings of $K_3$ and $K_4$ in the real projective plane, and the antipodal quotients of the icosahedron (on the sphere) and the great dodecahedron (on an orientable surface of genus $4$), giving regular embeddings of $K_6$ in the real projective plane and in a nonorientable surface of genus $5$.
\end{prop}

These nonorientable regular embeddings of $K_6$ form a Petrie dual pair, of types $\{3, 5\}_5$ and $\{5, 5\}_3$; they have automorphism group $L_2(5)\cong A_5$, and the second of them appears as entry N5.3 in Conder's computer-generated list of regular maps~\cite{C}.

We now prove the existence part of Theorem~\ref{main-theorem}. Case~(1) is easily dealt with, as follows.

\begin{lemma} \label{nonori-Hamming-n=2}
There is a nonorientable regular embeddings of $H(2,2)$ in the real projective plane. It has type $\{8, 2\}_8$, and its automorphism group is a dihedral group of order $16$.
\end{lemma}
\begin{pf}
Since $H(2,2)$ is a cycle of length $4$, we obtain a nonorientable regular embedding of this graph in the real projective plane by taking the antipodal quotient of  the regular embedding of a cycle of length $8$ in the sphere. The resulting map has one octagonal face, and its vertices have valency $2$, so it has type $\{8, 2\}_8$. Its automorphism group is the symmetry group of an octagon, namely a dihedral group of order $16$.
\end{pf}

For case~(2) we use the following result, which forms part of the first author's classification~\cite{J} of the orientably regular embeddings of $H(d,n)$.

\begin{prop} \label{orientable-Hamming}
If $d \ge 2$ and $n=3$ or $4$ there is an orientable regular embedding of $H(d,n)$. If $n=3$ it has type $\{m, 2d\}_6$, where the covalency $m$ is $2d$ or $3d$ as $d$ is even or odd. If $n=4$ it has type $\{3d, 3d\}_4$. Its automorphism group is a semidirect product of an elementary abelian normal subgroup of order $n^d$, acting regularly on the vertices, by a dihedral group of order $2d(n-1)$ fixing a vertex and acting naturally on its neighbours.
\end{prop}

(In fact, it is shown in~\cite{J} that these are the only orientable regular embeddings of Hamming graphs for $n\ge 3$.)

\begin{corollary} \label{nonorientable-Hamming}
If $d \ge 1$ and $n=3$ or $4$, there is a nonorientable regular embedding of $H(d,n)$. If $n=3$ it has genus $(2d-3)3^{d-1} + 2$, and it has type type $\{6, 2d\}_m$ where $m$ is $2d$ or $3d$ as $d$ is even or odd. If $n=4$ it has genus $(3d-4)4^{d-1} + 2$ and type $\{4, 3d\}_{3d}$. The automorphism group and its action on the vertices are as described in Proposition~\ref{orientable-Hamming}.
\end{corollary}
\begin{pf}
Let $d \ge 1$ and $n=3$ or $4$, and let $\cal M$ be the orientable regular embedding of $H(d,n)$ given by Proposition~\ref{orientable-Hamming}. Since Hamming graphs are not bipartite for $n\ge 3$, the Petrie dual ${\cal N}=P({\cal M})$ of $\cal M$ is a nonorientable regular embedding of $H(d,n)$. Petrie duality transposes covalency and Petrie length, so the type of $\cal N$ is as claimed, and the genus then follows from the Euler formula. Since Petrie duality preserves the automorphism group and its action on vertices, $\Aut{\cal N}$ is as in Proposition~\ref{orientable-Hamming}.
\end{pf}

For example, consider the first nontrivial case, namely $d=2$. If $n=3$ then $\cal N$ has genus $5$ and type $\{6,4\}_4$, and is the dual of entry N5.2 in Conder's list~\cite{C} of nonorientable regular maps; if $n=4$ then $\cal N$ has genus $10$ and type $\{4,6\}_6$, and is entry N10.1 in this list. If $d=3$ then for $n=3$ we obtain N29.2, of genus $29$ and type $\{6,6\}_9$, and for $n=4$ we have N82.1, of type $\{4,9\}_9$ and genus $82$.

The rest of this section is devoted to case~(3) of Theorem~\ref{main-theorem}, so we take $n=6$.

For $d=1$, the Hamming graph $H(1,6)$ is the complete graph $K_6$. It is known by work of James~\cite{Ja} and Wilson~\cite{Wi} that there are, up to isomorphism, two nonorientable regular
embeddings of $K_6$, described in Proposition~\ref{regKn}(b) and the remarks following it. Here we will give a group-theoretic
construction of a Petrie dual pair of nonorientable regular embeddings of $H(2,6)$.

Just as the nonorientable regular embeddings of $H(1,6)$ can be
obtained from the action of $A_5$ by conjugation on
its six Sylow 5-subgroups, those of $H(2,6)$ can be obtained from
the corresponding action of the group $P:=PGL_2(9)$. This has a
simple subgroup $L:=L_2(9)\cong A_6$ of index 2, which has twelve
icosahedral subgroups $A\cong A_5$, forming two conjugacy classes
$I$ and $J$ of size six: those in $I$ are the point stabilisers in
the natural action of $A_6$, while those in $J$ act transitively
as $L_2(5)$ in its natural action on the six points of the
projective line ${\mathbb P}^1(5)$ over the field ${\mathbb F}_5$. These two classes are
transposed by conjugation by elements of $P \setminus L$, so they
merge to form a single conjugacy class in $P$. Each of the 36
Sylow 5-subgroups $S \cong C_5$ of $P$ lies in exactly one group
$A \in I$ and one group $B \in J$, so by defining two Sylow
5-subgroups to be adjacent if they are contained in a common
icosahedral subgroup, we obtain a graph $H\cong H(2,6)$ on which $P$
acts by conjugation as a group of automorphisms. The stabiliser of
a vertex $S$ is its normaliser $N_P(S)$ in $P$, a dihedral group
$D$ of order 20 which acts transitively on the ten neighbours of
$S$. We will use this to construct nonorientable
regular embeddings of $H(2,6)$.

Let $\lambda, \rho$ and $\tau$ be the
elements of $P$ corresponding to the matrices
$$ M_{\lambda} =\left( \begin{array}{cc}
-1 &1 \\
1 & 1  \end{array} \right),~~ M_{\rho} =\left( \begin{array}{cc}
0 & 1+i \\
-1 & 0  \end{array} \right)~~\mbox{and}~~ M_{\tau} =\left( \begin{array}{cc}
1 &1 \\
1 & -1  \end{array} \right)$$
in $GL_2(9)$, where ${\mathbb F}_9 = {\mathbb F}_3(i)$ with $i^2=-1$. The elements
$\rho$ and $\tau$ generate a dihedral group $D$ of order 20, which
is maximal in $P$ and does not contain $\lambda$, so these three
elements generate $P$. They satisfy
$$\lambda^2=\rho^2=\tau^2=(\lambda\tau)^2=1,$$
so they determine a regular map $\cal M$ with $\Aut (\cal M) \cong
P$: the vertices, edges and faces correspond to the cosets in $P$
of the subgroups $\langle \rho , \tau \rangle$, $\langle \lambda , \tau
\rangle$ and $\langle \lambda , \rho \rangle$, with incidence given by
non-empty intersection. Since $\lambda\rho$, $\rho\tau$ and $\lambda\rho\tau$
have orders $10, 10$ and $8$ respectively, $\cal M$ has type $\{
10, 10\}_8$, while its Petrie dual $\cal N$ has type $\{ 8, 10
\}_{10}$. Since $|P|=720$, $\cal M$ and $\cal N$ have $36$ and $45$ faces, respectively, so they have Euler characteristics $\chi=-108$ and $-99$. Since $\det M_{\lambda}=1$, $\lambda$ is contained in the unique
subgroup $L$ of index 2 in $P$, so $\cal M$ and $\cal N$ are
nonorientable and hence have genera $2-\chi=110$ and $101$ (they appear as entries N110.7 and N101.8 in~\cite{C}). These two maps have the same underlying graph $K$, which we will show is isomorphic to $H$.

The subgroup $D=\langle \rho, \tau \rangle$ stabilising a vertex of
$K$ is the normaliser in $P$ of the Sylow 5-subgroup $S =
\langle a  \rangle $, where $a= (\rho\tau)^2$; thus the vertices of
$K$ can be identified with the 36 Sylow 5-subgroups of $P$,
permuted by conjugation, and hence with the vertices of $H$. Now
we consider the edges of $K$. The Sylow 5-subgroup $S$
corresponds to the vertex $D$ of $K$, and this is adjacent in
$K$, through their common edge $\langle \lambda, \tau \rangle$, to
the vertex corresponding to the Sylow 5-subgroup $T=S^{\lambda}$
generated by the element $b=a^{\lambda}$. The subgroup $A$ generated
by $a$ and $b$ is also generated by $a$ and $a^2b$, which have
orders $5$ and $3$, while their product has order $2$; thus $A$ is
an epimorphic image of the triangle group of type $(5,3,2)$,
isomorphic to the simple group $A_5$, so $S$ and $T$ generate an
icosahedral subgroup of $P$. Thus $S$ and $T$ are neighbours in
$H$, and since $P$ acts edge transitively on both $K$ and
$H$, it follows that their edge sets correspond. Thus our
identification of their vertex sets gives an isomorphism $K
\rightarrow H$ commuting with the actions of $P$ on these two
graphs. Since $\cal M$ and $\cal N$ were constructed as nonorientable
regular embeddings of $K$, this shows that they also yield
such embeddings of $H\cong H(2,6)$.

To summarise, we have proved the following result.

\begin{lemma} \label{nonori-Hamming-n=6}
There is a Petrie dual pair of nonisomorphic nonorientable regular embeddings of $H(2,6)$. They have genera $110$ and $101$, and types $\{10,10\}_8$ and $\{8,10\}_{10}$. Their automorphism group is isomorphic to $PGL_2(9)$.
\end{lemma}

In the next section, we will show that the regular maps
constructed in this section are, up to isomorphism, the only nonorientable regular
embeddings of Hamming graphs.

\section{Classification of nonorientable Hamming maps}

In this section, we classify the nonorientable regular embeddings of
Hamming graphs $H(d,n)$  by completing the proof of Theorem \ref{main-theorem}. The case $d=1$, where $H(1,n)\cong K_n$, is covered by the work of James~\cite{Ja} and Wilson~\cite{Wi} summarised in Proposition~\ref{regKn}(b), showing that $n=3, 4$ or $6$, so from now on we will assume that $d\ge 2$. If $n=2$ then $H(d,n)$ is isomorphic to the $d$-cube $Q_d$, and it is known from work of Nedela and the second author~\cite{KN} that the only nonorientable regular embedding of $Q_d$ for any $d\ge 2$ is the embedding of $Q_2 \cong H(2,2)$ described in Lemma~\ref{nonori-Hamming-n=2}. This deals with the case $n=2$, including part~(1) of Theorem~\ref{main-theorem}, so we will assume from now on that $n \ge 3$.

For notational convenience we will regard $[n]=\{0, \ldots, n-1\}$ as the ring ${\mathbb Z}_n$, and we will label the coordinate places of vertices with elements $i\in[d]:=\{0, 1, \ldots, d-1\}$. We will identify $\Aut(H(d,n))$ with $S_n \wr S_d$, acting as described at the beginning of Section~2, and we will write each element of this group in the form $\delta = \delta^*\delta'$ where $\delta^*=(\delta_0,\ldots, \delta_{d-1})\in S_n^d$ with each $\delta_i\in S_n$, and $\delta'\in S_d$. We will use the facts that the projection $\delta \mapsto \delta'$ is a homomorphism, and that the elements $\delta$ with each $\delta_i$ even form a subgroup $A_n\wr S_d$ of $S_n\wr S_d$.

Let us define permutations
$$\alpha_d = (0 ~ 1~ \ldots ~ d-1) \quad {\rm and} \quad
\beta_d = (0) (1~  d-1)(2~d-2) \ldots $$
in $S_d$ (so that $\beta_2$ is the identity permutation $id$), and
$$\beta_n = (0) (1~  n-1)(2~n-2) \ldots  \quad {\rm and} \quad \gamma_n = (1 ~ 2~ \ldots ~ n-1)$$
in $S_n$. Let $O$ denote the vertex $(0, \ldots, 0)$, and for each $i=0, \ldots, d-1$ let $e_i$ denote the vertex $(0,\ldots, 0,1,0,\ldots, 0)$ having a single non-zero coordinate $1$ in position $i$. By using the natural structure of the vertex set ${\mathbb Z}_n^d$ as a ${\mathbb Z}_n$-module, we can denote a general vertex $(v_0, v_1,\ldots, v_{d-1})$ of $H(d,n)$ by $\sum_{i=0}^{d-1}v_i  e_i$.

\begin{lemma} \label{form-admissible-conjugacy}
For any admissible triple $(\lambda , \rho ,\tau)$, corresponding to an orientable or nonorientable regular embedding of $H(d,n)$ with $d\ge 2$ and $n\ge 3$, there exists an automorphism $\phi$ of $H(d,n)$ such that
\begin{eqnarray*}
\phi^{-1}\tau\phi  &=&  ((0)(1)(2~~n-1)(3~~n-2)\ldots, \beta_n, \ldots, \beta_n)\beta_d \\
\phi^{-1}\rho\tau\phi &=& (id,id, \ldots, \gamma_n)\alpha_d \\
\phi^{-1}\lambda\tau\phi &=& (\sigma_0 ,\sigma_1, \ldots, \sigma_{d-1})\theta,
 \end{eqnarray*}
 where $0^\theta = 0$, $\theta^2 = id$, $0^{\sigma_0} =1$, $1^{\sigma_0} =0$, and $0^{\sigma_i} =0$ and $\sigma_i \sigma_{i^\theta} = id$ for each $i \in [d] \setminus \{0 \}$. In particular, if $(\lambda , \rho ,\tau)$ is a nonorientable admissible triple for
$H(d,n)$, we can assume that $\theta = \beta_d$ in the above form. The subgroup $\langle (\rho\tau)^d, \lambda\tau, \tau \rangle$ is isomorphic to the automorphism group of a regular embedding of $K_n$, which is nonorientable  if $d \geq 3$.
 \end{lemma}
 
\begin{pf}
Let  $(\lambda , \rho ,\tau)$ be an admissible triple for $H(d,n)$, generating the automophism group $G$ of a regular embedding $\cal M$ of $H(d,n)$. Thus $\rho$ and $\tau$ fix a vertex $v$, which $\lambda$ transposes with an adjacent vertex $w$ fixed by $\tau$. Since $\Aut(H(d,n)) = S_n \wr S_d$ acts transitively on the arcs of $H(d,n)$, there is a graph automorphism $\phi$ sending $v$ and $w$ to $O$ and $e_0$. Then $\lambda_1 := \phi^{-1}\lambda\phi$ transposes $O$ and $e_0$, while $\rho_1 := \phi^{-1}\rho\phi$ fixes $O$, and $\tau_1 := \phi^{-1}\tau\phi$ fixes $O$  and $e_0$. The conjugate admissible triple $(\lambda_1, \rho_1, \tau_1)$ generates the automorphism group $G_1:=\phi^{-1}G\phi$ of a regular embedding ${\cal M}_1=\phi({\cal M})$ of $H(d,n)$. Now $\langle \rho\tau\rangle$ acts transitively on the $d(n-1)$ neighbours of $v$; these form $d$ disjoint cliques, which are blocks of imprimitivity for this group, so $(\rho\tau)'$ is a cyclic permutation of $[d]$, and we may choose $\phi$ so that $(\rho_1\tau_1)'=\alpha_d$.  More specifically, we may choose $\phi$ (or equivalently relabel the vertices) so that $\lambda_1$, $\rho_1$ and $\tau_1$ satisfy 
\begin{enumerate}
\item [{\rm (i)}]
$\rho_1\tau_1 = (id,id, \ldots, \gamma_n)\alpha_d$,
so that for any $k \in {\mathbb Z}_n \setminus\{0\}$ we have $(k e_i)^{\rho_1\tau_1} = ke_{i+1}$ for all $i= 0,\ldots, d-2$, while
$(ke_{d-1})^{\rho_1\tau_1}$ is $(k+1)e_0$ or $e_0$ as $k \neq n-1$ or $k=n-1$ respectively,
\item [{\rm (ii)}] $\lambda_1$ transposes the vertices $O$ and $e_0$, and preserves the face incident with $O$, $e_0$ and $e_{n-1}$,
\item [{\rm (iii)}] $\rho_1$ fixes the vertex $O$ and preserves the face incident with $O$, $e_0$ and $e_{n-1}$, and
\item [{\rm (iv)}] $\tau_1$ fixes both $O$  and $e_0$.
\end{enumerate}
Since $\tau_1$ inverts $\rho_1\tau_1$ by conjugation, and commutes with $\lambda_1$, one can easily verify that
\begin{eqnarray*}
\tau_1 &=&  ((0)(1)(2~~n-1)(3~~n-2)\cdots,\beta_n, \ldots, \beta_n)\beta_d \\
\lambda_1\tau_1 &=& (\sigma_0 ,\sigma_1, \ldots, \sigma_{d-1})\theta,
 \end{eqnarray*}
 where $0^\theta = 0$, $\theta^2 = id$, $0^{\sigma_0} =1$, $1^{\sigma_0} =0$, and $0^{\sigma_i} =0$ and $\sigma_i \sigma_{i^\theta} = id$ for all $i \in [d] \setminus \{0 \}$. Note that if $d=2$ then both $\beta_d$ and $\theta$ are the identity permutation.
 
Let $R_1 :=\rho_1 \tau_1$ and $L_1 := \lambda_1\tau_1$. The
subgraph $K$ of $H(d,n)$ induced by the set of vertices $\{ke_0 ~|~ k \in {\mathbb Z}_n\}$ is isomorphic to the complete graph $K_n$. It is invariant under $R_1^d$, $L_1$ and $\tau_1$, and hence under the subgroup  $H= \langle R_1^d, L_1, \tau_1 \rangle$ of $G_1$ which they generate. Now $H$ acts transitively (and hence regularly) on the flags of ${\cal M}_1$ incident with $K$, with $(\lambda_1, R_1^d\tau_1, \tau_1)$ an admissible triple, so $H$ is isomorphic to the automorphism group of a regular embedding $\cal K$ of $K_n$. (In fact, we have constructed $\cal K$ from ${\cal M}_1$ by applying Wilson's map operation $H_d$~\cite{Wi79}, which raises the local rotation of arcs around each vertex to its $d$-th power; in this case $d$ is not coprime to the valency $d(n-1)$, so $H$ is not transitive on the flags of ${\cal M}_1$, and we have taken a single orbit to define $\cal K$.)

Since $n \ge 3$ the graph $K$ contains a cycle of length $3$, so there exist $k_1, k_2 , k_3 \in {\mathbb Z}_{n-1} \setminus \{0 \}$ such that the element
$$\delta := L_1R_1^{dk_1}L_1R_1^{dk_2}L_1R_1^{dk_3}$$
of $H$ is in $\langle \tau_1 \rangle$. Now
$$R_1^d = (\gamma_n, \ldots, \gamma_n)\in S_n\times\cdots\times S_n,$$
so $\delta$ has image $\delta'=\theta^3=\theta$ in $S_d$.

We will assume that $(\lambda , \rho ,\tau)$ is a nonorientable admissible triple, so the triple $(\lambda_1 , \rho_1 ,\tau_1)$ is also nonorientable. Thus $\tau_1 \in  \langle R_1, L_1 \rangle$, so by taking images in $S_d$ we have $\tau_1'\in\langle R_1', L_1'\rangle$, that is, $\beta_d \in \langle\alpha_d, \theta\rangle$. If $d \ge 3$ then visibly $\beta_d \not\in \langle\alpha_d\rangle$, so $\theta \neq id$.  Thus $\delta'\ne id$, so $\delta\ne id$ and hence  $\delta = \tau_1$  (and $\theta = \tau_1' = \beta_d$). Thus $\tau_1$ is contained in the subgroup $\langle R_1^d , L_1 \rangle$ of $H$, so $\langle R_1^d , L_1 \rangle = H$ and $\cal K$ is a nonorientable regular embedding of $K_n$.
\end{pf}

It is thus sufficient to consider nonorientable admissible triples $(\lambda, \rho ,\tau)$ for $H(d,n)$ of the form given by $(\lambda_1, \rho_1, \tau_1)=(\phi^{-1}\lambda\phi, \phi^{-1}\rho\phi, \phi^{-1}\tau\phi)$ in Lemma~\ref{form-admissible-conjugacy}. For any such triple $(\lambda, \rho ,\tau)$, the group $H=\langle (\rho\tau)^d , \lambda\tau, \tau \rangle$ is isomorphic to the automorphism group of a regular embedding $\cal K$ of $K_n$, which is nonorientable if $d\ge 3$, but may be orientable or nonorientable if $d=2$. In either case, since we are assuming that $n\ge 3$, Proposition~\ref{regKn} implies the following:

\begin{corollary}\label{cor:n=346}
If there is a nonorientable regular embedding of $H(d,n)$ with $d\ge 2$ and $n\ge 3$, then $n=3, 4$ or $6$. 
\end{corollary} 

To deal with case~(2) of Theorem~\ref{main-theorem}, let us assume that $n=3$ or $4$.

 \begin{lemma} \label{unique-nonori-Hamming-n=3,4}
If $d \ge 2$ and $n=3$ or $4$, there is, up to isomorphism, at most one nonorientable regular embedding of $H(d,n)$.
\end{lemma}

\begin{pf}
Let $n=3$ and $d=2$. Then $\tau = (id,(1~~2))id$, $R=(id,(1~~2))(0~~1)$ and $L = ((0~~1), \sigma_1 )id$, where $\sigma_1 =id$ or $(1~~2)$. One can check that   $LR^2LR^2LR^2 = \tau$ or $id$ as $\sigma_1 = id$ or $(1~~2)$ respectively. By Proposition~\ref{orientable-Hamming} there is an orientable regular embedding of $H(2,3)$, which must correspond to the case $\sigma_1 =(1~~2)$, so there is at most one nonorientable regular embedding of $H(2,3)$, with $\sigma_1=id$.

Now let $n=3$ and $d \geq 3$.  Then $\sigma_0 = (0~~ 1)$, and for each $i=1, 2, \ldots, d-1$ the element $\sigma_i = \sigma_{d-i}$ is the identity or $(1~~2)$ because $0^{\sigma_i} =0$. By Lemma~\ref{form-admissible-conjugacy}, the subgroup $\langle (\rho\tau)^d, \lambda\tau \rangle$ is isomorphic to the automorphism group of a nonorientable regular embedding $\cal K$ of $K_3$. By Proposition~\ref{regKn}(b) the underlying surface of $\cal K$ is the real projective plane, and the graph $K_3$ has a neighbourhood homeomorphic to a M\"obius band. Thus $LR^dLR^dLR^d = \tau$, which implies that
 $$\sigma_i (1~~2)\sigma_{d-i}(1 ~~2)\sigma_i (1 ~~2) = \beta_3  =(1 ~~2)$$
 for each $i=1,2,\ldots, d-1$, and hence $\sigma_i (1~~2)\sigma_{d-i}(1 ~~2)\sigma_i =id$. This implies that  $\sigma_i = \sigma_{d-i} = id$, so there is at most one nonorientable regular embedding of $H(d,3)$ for each $d \ge 3$.

Next let $n=4$ and $d=2$. Then $\tau = ((2~~3),(1~~3))id$, $R=(id,(1~~2~~3))(0~~1)$ and $L = (\sigma_0, \sigma_1 )id$. Now $\tau\in\langle R, L \rangle$; since the coordinates of the element $\tau^*\in S_4\times S_4$ are both odd permutations, whereas those of $R^*$ are both even, each $\sigma_i$ must be odd. Since $0^{\sigma_1} =0$, it follows that $\sigma_1$ is a transposition.  By Lemma~\ref{form-admissible-conjugacy} the subgroup $H=\langle R^2 , L, \tau \rangle$ is isomorphic to the automorphism group of a regular embedding $\cal K$ of $K_4$. Now $K_4$ contains a 3-cycle, so there exist $k_1, k_2 , k_3 \in {\mathbb Z}_3 \setminus \{0 \}$ such that the element
$$\delta := LR^{2k_1}LR^{2k_2}LR^{2k_3}$$
is in $\langle \tau \rangle$. Since the second coordinate of $\delta^*$ is an odd permutation, $\delta = \tau$, so $\tau \in \langle R^2, L \rangle$ and $\cal K$ is nonorientable. Thus $\cal K$ is the unique nonorientable regular embedding of $K_4$, namely the antipodal quotient of the cube. This implies that $\sigma_0 = (0~~1)$ and  $LR^2LR^4LR^2 =
 \tau$. Comparing the second coordinates of $(LR^2LR^4LR^2)^*$ and $\tau^*$, we have $\sigma_1 (1~~2~~3)\sigma_1(1 ~~3~~2)\sigma_1(1 ~~2~~3) = (1
 ~~3)$. Among the three possible transpositions $\sigma_1\in S_4$ fixing $0$, one easily checks that only $(1~~3)$ satisfies this equation, so there is at most one nonorientable regular embedding of $H(2,4)$.

Finally let $n=4$ and $d \geq 3$. By Lemma~\ref{form-admissible-conjugacy},  $\langle R^d , L \rangle$ is isomorphic to the automorphism group of a nonorientable regular embedding of $K_4$, which again implies that $\sigma_0 = (0~~1)$ and $LR^dLR^{2d}LR^d = \tau$. Hence for each $i=1,2,\ldots, d-1$,
$$\sigma_i (1~~2~~3)\sigma_{d-i}(1 ~~3~~2)\sigma_i(1 ~~2~~3)
= \beta_4  =(1 ~~3).$$
Since  $\sigma_i \sigma_{d-i} =id$, both $\sigma_i$ and $\sigma_{d-i}$ are odd permutations. Since they fix $0$, this implies that both $\sigma_i$ and $\sigma_{d-i}$ are transpositions, and hence $\sigma_i = \sigma_{d-i}$.
Now one can easily show that $\sigma_i = \sigma_{d-i}=(1~~3)$, so there is at most one nonorientable regular embedding of $H(d,4)$ for each $d \ge 3$.
\end{pf}

This result and Corollary~\ref{nonorientable-Hamming} deal with the cases $n=3$ and $n=4$, including part~(2) of Theorem~\ref{main-theorem}. By Corollary~\ref{cor:n=346}, in order to complete the proof, and to deal with part~(3), we may assume from now on that $n=6$.

\begin{lemma} \label{two-nonori-Hamming-n=6}
There is no nonorientable regular embedding of $H(d,6)$ for any $d\ge 3$, and there are, up to isomorphism, at most two nonorientable regular embeddings of $H(2,6)$.
\end{lemma}
\begin{pf}
We will use the easily verified fact that the two nonorientable regular embeddings of $H(1,6)\cong K_6$ described in Proposition~\ref{regKn}(b) are derived from admissible triples $(\lambda, \rho ,\tau)$ for this graph with
$$\tau = (2~~5)(3~~4),~ ~~
R = \gamma_6 = (1~~2~~3~~4~~5) ~~~\mbox{and}~~~
L=(0~~1)(2~~5) ~~ \mbox{or}~~ (0~~1)(3~~4)$$
in $S_6$. If $L=(0~~1)(2~~5)$ then $LRLR^4LR^2LR^2 = \tau$, which implies that there is a 4-cycle in $K_6$ with a neighbourhood homeomorphic to a M\"obius band.  Since these two nonorientable embeddings form a Petrie dual pair, it follows that both embeddings have such a 4-cycle.

Now let  $(\lambda, \rho, \tau)$ be a nonorientable admissible triple for $H(d,6)$, with $d\ge 2$. By Lemma~\ref{form-admissible-conjugacy}, we can assume that
\begin{eqnarray*}
\tau  &=&  ((2~~5)(3~~4)\cdots,\beta_6, \ldots, \beta_6)\beta_d \\
R  &=& \rho\tau = (id,id, \ldots, \gamma_6)\alpha_d \\
 L  &=& \lambda\tau = (\sigma_0 ,\sigma_1, \ldots, \sigma_{d-1})\beta_d.
 \end{eqnarray*}
Moreover the subgroup $H = \langle R^d , L, \tau \rangle$ is isomorphic to the automorphism group of a regular embedding $\cal K$ of the subgraph $K\cong K_6$ induced by the vertices $ke_0$ for $k\in{\mathbb Z}_6$. Proposition~\ref{regKn} shows that $\cal K$ must be nonorientable, so by our earlier argument there is a 4-cycle in $K$ with a neighbourhood in $\cal K$ homeomorphic to a M\"obius band. Thus there exist $k_1, k_2, k_3$ and $k_4$ such that $LR^{dk_1} LR^{dk_2} LR^{dk_3} LR^{dk_4}=\tau$.  By taking images in $S_d$ we see that $\beta_d = id$, and hence $d = 2$. Thus there is no nonorientable regular embedding of $H(d,6)$ for any $d \ge 3$.

Finally, let $d=2$.  By considering the actions of $R^2=(\gamma_6, \gamma_6)$ and $L$ on the coefficients $k\in{\mathbb Z}_6$ of the vertices $ke_0$ of $K$ we see that $\gamma_6$ and $\sigma_0$ generate a subgroup of $S_6$ isomorphic to $\Aut{\cal K}$. The two regular embeddings of $K_6$ in Proposition~\ref{regKn}(b) have automorphism groups isomorphic to $L_2(5)$, so $\langle \gamma_6,\sigma_0\rangle\cong L_2(5)$. This group is perfect, so $\sigma_0$ is an even permutation; since $\sigma_0$ is an involution, it is a double transposition. Now $\sigma_0$ transposes $0$ and $1$, so there are just six possibilities, and inspection shows that $\langle\gamma_6, \sigma_0\rangle \cong L_2(5)$ if and only if $\sigma_0 = (0~~1)(2~~5)$ or $(0~~1)(3~~4)$. 

If $\sigma_0 = (0~~1)(2~~5)$ then $(LR^2)^3 =id$, so $\cal K$ has covalency $3$; by Proposition~\ref{regKn}(b) it is therefore the antipodal quotient of the icosahedral map. This satisfies $LR^4LR^6LR^4 = \tau$, and these two equations imply that $\sigma_1$ must satisfy
$$(\sigma_1(1~~2~~3~~4~~5))^3 =id$$
and
$$\sigma_1(1~~3~~5~~2~~4)\sigma_1(1~~4~~2~~5~~3)\sigma_1(1~~3~~5~~2~~4)=(2~~5)(3~~4).$$
By the first equation, $\sigma_1$ is an even permutation. Note also that $\sigma_1$ is an involution fixing $0$. One can check that $(1~~2)(4~~5)$ is the only permutation in $A_6$ satisfying these conditions on $\sigma_1$, so there is at most one nonorientable regular embedding of $H(2,6)$ with $\sigma_0 = (0~~1)(2~~5)$.

If $\sigma_0 = (0~~1)(3~~4)$ then a similar argument, in which $\cal K$ is now the antipodal quotient of the great dodecahedron, shows that
the only possibility for $\sigma_1$ is $(1~~4)(2~~5)$, so there is at most one nonorientable regular embedding of $H(2,6)$ with $\sigma_0 = (0~~1)(3~~4)$.
 \end{pf}

In the case $d=2$ of the above proof, if $\sigma_0 = (0~~1)(2~~5)$ and $\sigma_1 = (1~~2)(4~~5)$ the covalency of the corresponding embedding of $H(2,6)$ is $8$, and if $\sigma_0 = (0~~1)(3~~4)$ and $\sigma_1 = (1~~4)(2~~5)$ it is $10$, so these embeddings are respectively the maps $\cal N$ and $\cal M$ of types $\{8, 10\}_{10}$ and $\{10, 10\}_8$ constructed in the proof of Lemma~\ref{nonori-Hamming-n=6}.

This deals with part~(3) of Theorem~\ref{main-theorem}, and completes the proof of this result.

 \end{document}